\documentclass{amsart}

\usepackage{amsmath}
\usepackage{amssymb}
\usepackage{amsthm}
\usepackage[matrix,arrow,graph,curve,frame]{xy}
\usepackage{hyperref}
\usepackage{mathrsfs}
\usepackage{enumerate}

\theoremstyle{plain}      \newtheorem{theorem}{Theorem}
\theoremstyle{plain}      
\theoremstyle{definition} \newtheorem{example}[theorem]{Example}
\theoremstyle{remark}     \newtheorem*{remark}{Remark}
\theoremstyle{remark}     \newtheorem*{note}{Note}

\newcommand{\A}{\ensuremath{\mathscr{A}}}

\newcommand{\V}{\ensuremath{\mathscr{V}}}
\newcommand{\R}{\ensuremath{\mathbb{R}}}
\newcommand{\G}{\ensuremath{\mathcal{G}}}
\newcommand{\ox}{\ensuremath{\otimes}}
\newcommand{\op}{\ensuremath{\mathrm{op}}}
\newcommand{\ob}{\ensuremath{\mathrm{ob}}}
\newcommand{\olstar}{\ensuremath{\mathbin{\overline{*}}}}
\newcommand{\ulstar}{\ensuremath{\mathbin{\underline{*}}}}
\newcommand{\x}{\ensuremath{\times}}

\newcommand{\Vect}{\ensuremath{\mathbf{Vect}}}

\newcommand{\ra}{\ensuremath{\xymatrix@C=4.5ex@1{\ar[r]&}}}
\newcommand{\Ra}{\ensuremath{\xymatrix@C=3ex@1{\ar@{=>}[r]&}}}

\title{$*$-Autonomous categories in quantum theory}
\author{Brian Day}
\address{Centre of Australian Category Theory, Macquarie University, NSW,
2109, Australia}
\date{November 5, 2011}

\begin{document}
\maketitle

\begin{abstract}
$*$-Autonomous categories were initially defined by M.\ Barr to describe a type of duality carried by many monoidal closed categories. Later they were generalised by the current author to include $*$-autonomous promonoidal categories. Together, these structures under ``convolution'' product give a clear indication of the usefulness of $*$-autonomy in quantum mathematics and related areas.
\end{abstract}

\section{Introduction}

So-called ``$*$-autonomous'', or ``Frobenius'', category structures occur widely in
mathematical quantum theory. This trend was observed in~\cite{4}, mainly in
relation to Hopf algebroids, and continued in~\cite{9} with a general account
of Frobenius monoids.

Below we list some of the $*$-autonomous partially ordered sets $\A =
(\A,p,j,S)$ that appear in the literature, an abstract definition of
$*$-autonomous promonoidal structure being made in~\cite[\S 7]{4}. Without
going into much detail, we also note some features of the convolution $[\A,\V]$
(defined in~\cite{2}) of a given such $\A$ with a complete $*$-autonomous
monoidal category $\V$. A monoidal functor category of this type is a
completion of $\A$, with an appropriate universal property; it is always again
$*$-autonomous (as seen in~\cite{4} for example).

The basic descriptions of promonoidal (equals premonoidal) structure and the
resulting convolution product are given in~\cite{2} and~\cite{4}. We shall
more-or-less regard the less complete promonoidal structures as the
microscopic versions of their convolution counterparts.

The examples discussed below are mostly based in the extended positive real
numbers $\R^\infty_{\ge 0}$ with the $*$-autonomous monoidal structure of
multiplication and identity 1 (we simply define $\infty \ox 0 = 0$, and
$\infty \ox r = \infty$ for $r \neq 0$).

\begin{remark}
The process of adding $\infty$ to $\R_{\ge 0}$ is quite general. For example,
one can add 0 (initial) and $\infty$ (terminal) to any partially ordered group
$G$ and obtain $G_0^\infty$ which is a $*$-autonomous monoidal category.
Furthermore, the group $G$ can be replaced by any rigid closed category, etc.
\end{remark}

The category $\R^\infty_{\ge 0}$ is isomorphic, under exponentiation, to the
$*$-autonomous category
\[
    \R^\infty_{-\infty} = \{-\infty\} \cup \R \cup \{\infty\}
\]
with the monoidal structure of addition and identity 0. In the following, each
poset $\A$ is viewed as a category under
\[
    \A(a,b) =
    \begin{cases}
    1 & \text{iff } a \le b \\
    0 & \text{else,}
    \end{cases}
\]
and the base category is $\V = \R^\infty_{\ge 0}$ unless otherwise mentioned.

\section{Examples}

\begin{example}[Submodular functions~\cite{8}]
Let $E$ be a set and $\A = \mathcal{P}(E)$ (discrete). For $\V = \R_{\ge
0}^\infty$, let
\[
    p(a,b,c) =
    \begin{cases}
    1 & \text{iff } (a \cup b = c \text{ and } a \cap b = \emptyset)
    \text{ iff } a+b=c \\
    0 & \text{else,}
    \end{cases}
\]
and
\[
    j(a) =
    \begin{cases}
    1 & \text{iff } a = \emptyset \\
    0 & \text{else.}
    \end{cases}
\]
Then $(\A,p,j,S)$ becomes a $*$-autonomous promonoidal $\V$-category if we take
\[
    Sa = E -a,
\]
since $p(a,b,Sc) = 1$ iff $a+b=E-c$, i.e., $a+b+c = E$.
\end{example}

\begin{note}
In fact Narayanan uses $\V = (\R^\infty_{-\infty},+,0)$ instead of
$\R_{\ge 0}^\infty$, and discusses the \emph{upper convolution} on $[\A,\V]$
given by
\begin{align*}
f \olstar g(c) &= \sup_{ab} f(a)+g(b) + p(a,b,c) \\
               &= \sup_{a \subset c} f(a)+g(c-a).
\end{align*}

(since $\A = \mathcal{P}(E)$ is discrete, $[\A,\V]$ equals all functions from
$\mathcal{P}(E)$ to $\R_{-\infty}^\infty$), and also the \emph{lower
convolution}
\begin{align*}
f \ulstar g(c) &= (\sup_{ab} f(a)^* + g(b)^* + p(a,b,c))^* \\
               &= -(\sup_{ab} -f(a) - g(b) + p(a,b,c)) \\
               &= \inf_{a \subset c} f(a)+g(c-a).
\end{align*}
Of course in both cases we use the fact that
\[
    p(a,b,c) =
        \begin{cases}
        0 & \text{iff } a+b=c \\
        -\infty & \text{else}
        \end{cases}
\]
in $\R_{-\infty}^\infty$.

Note also that $f \in [\A,\V]$ is an upper convolution monoid iff
\[
    f(a+b) \ge f(a)+f(b) \text{ and } f(0) \ge 0,
\]
while $f$ is a lower convolution monoid iff
\[
    f(a+b) \le f(a)+f(b) \text{ and } f(0) \le 0.
\]
\end{note}

\begin{example}
Let $\A = (\A,\vee,\wedge,0,1,S)$ be an orthomodular lattice (see
Kalmbach~\cite{6} for example). Then the definitions
\[
    p(a,b,c) =
        \begin{cases}
        1 & \text{iff } a \vee b \le c \\
        0 & \text{else,}
        \end{cases}
\]
\[
    \A(a,b) =
        \begin{cases}
        1 & \text{iff } a \le b \\
        0 & \text{else,}
\end{cases}
\]
and
\[
    j(a) =
        \begin{cases}
        1 & \text{iff } a=0 \\
        0 & \text{else,}
        \end{cases}
\]
yield a $*$-autonomous promonoidal poset $(\A,p,j,S)$ where $Sa$ is the
orthocomplement of $a$ in $\A$. The orthomodularity law on $\A$ is equivalent
to either of the cyclic relations
\[
    p(a,b,Sc) = p(b,c,Sa) = p(c,a,Sb).
\]
\end{example}

\begin{example}[Browerian logics (Lawvere)]
Let $(\A, \le, \Ra, \wedge, 0, 1)$ be a Browerian logic. Define an $\R_{\ge
0}^\infty$-category by
\[
    \A(a,b) =
        \begin{cases}
        1 & \text{iff } a \le b \\
        0 & \text{else,}
        \end{cases}
\]
and let
\[
    p(a,b,c) =
        \begin{cases}
        1 & \text{iff } a \wedge b \le c \\
        0 & \text{else,}
        \end{cases}
\]
and
\[
    j(a) =
        \begin{cases}
        1 & \text{iff } a=1 \\
        0 & \text{else.}
        \end{cases}
\]
If we put $Sa=(a \Ra 0)$, then
\[
    p(a,b,Sc) = 1 \text{ iff } a \wedge b \wedge c = 0,
\]
hence $(\A,p,j,S)$ is a $*$-autonomous promonoidal category.
\end{example}

\begin{example}[Groupoids]
Let $\G$ be a groupoid and let $\A$ denote the set of arrows of $\G$. Define
\[
    p(a,b,c) =
        \begin{cases}
        1 & \text{iff } a b = c \\
        0 & \text{else,}
        \end{cases}
\]
\[
    j(a) =
        \begin{cases}
    1 & \text{iff } a \text{ is an identity} \\
    0 & \text{else,}
\end{cases}
\]
and $Sa = a^{-1}$. Then $p(a,b,Sc)=1$ iff $abc=1$, so $(\A,p,j,S)$ is
$*$-autonomous.
\end{example}

\begin{example}[(Non-commutative) probabilistic geometry (of~\cite{3})]
Let $\A$ be a poset with an associative promultiplication
\[
    p : \A^\op \x \A^\op \x \A \ra [0,1] \subset \R_{\ge 0}^\infty,
\]
where we interpret the value $p(a,b,c)$ as the probability that the point $c$
lies in the line joining $a$ and $b$. The convolution of poset maps $f$ and $g$
from $\A$ to $[0,1]$ is then the \emph{join} of $f$ to $g$:
\[
    f * g = \sup_{ab} f(a) g(b) p(a,b,-),
\]
while $f$ is \emph{convex} iff $f*f \le f$; i.e., iff $f$ is a convolution
semigroup.

In particular, note that if $\A$ is discrete, this $(\A,p,j)$ is $*$-autonomous
with respect to $Sa=a$ iff
\[
    p(a,b,c)=p(b,c,a)=p(c,a,b),
\]
and these simultaneously take the value $1$ iff the points $a$, $b$, and $c$
coincide.
\end{example}

\begin{example}[Generalized effect and difference algebras (cf.
Kalmbach~\cite{6}, Chapter 21; see also~\cite{5})]
Suppose the poset $\A$ has the structure of a (non-commutative, say)
generalized effect algebra $(\A,\oplus,0,\le)$. Let
\[
    p(a,b,c) =
        \begin{cases}
        1 & \text{iff } a \oplus b \le c \\
        0 & \text{else,}
    \end{cases}
\]
\[
    \A(a,b) =
        \begin{cases}
        1 & \text{iff } a \le b \\
        0 & \text{else,}
        \end{cases}
\]
and
\[
    j(a) =
        \begin{cases}
        1 & \text{iff } 0 \le a \\
        0 & \text{else,}
        \end{cases}
\]
then $(\A,p,j)$ is an (associative and unital) promonoidal category with some
extra properties (e.g., cancellation).

Similarly, if the poset $\A$ is a generalized \emph{commutative} difference
algebra $(\A,\ominus,0,\le)$ then
\[
    p(a,b,c) =
        \begin{cases}
        1 & \text{iff } a \le c \ominus b \\
        0 & \text{else,}
        \end{cases}
\]
and
\[
    j(a) =
        \begin{cases}
        1 & \text{iff } 0 \le a \\
        0 & \text{else,}
        \end{cases}
\]
yield an (associative and unital) promonoidal category $(\A,p,j)$. Note that a
\emph{commutative} generalized effect algebra is related to a (commutative)
generalized difference algebra by 
\[
    a \oplus b \le c \text{ iff } a \le c \ominus b.
\]

The key feature regarding~\cite{6} (Riecanov\'a) is that the promonoidal category
\[
    \mathscr{P} = \A + \A^\op
\]
constructed in the embedding theorem~\cite[Proposition 21.2.4]{6} (due to
Hedlikov\'a and Pulmannov\'a) is in fact $*$-autonomous under the switch map
\[
    S:(\A + \A^\op)^\op = \A^\op+\A \cong \A + \A^\op.
\]
Again the convolution
\begin{align*}
    [\mathscr{P},\V] &= [\A + \A^\op,\V] \\
                    &= [\A,\V] \x [\A^\op,\V]
\end{align*}
is $*$-autonomous monoidal and complete.
\end{example}

\begin{remark}
The construction of this $\mathscr{P}$ is closely related to (but seems not
the same as) that of the free $*$-autonomous promonoidal category on a given
promonoidal category due to Luigi Santocanale (unpublished?), the latter
giving the ``Chu construction'' upon convolution with $\V$.
\end{remark}

\begin{example}[Conformal field theory]
An example of a different nature $(\V=\Vect)$ arises in RCFT~\cite[4.17]{7}
as a $\Vect$ promonoidal structure on a (discrete) finite set $\A$ with a
distinguished base point 0.

The promultiplication
\[
    p:\A \x \A \x \A \ra \Vect_{\text{fd}}
\]
is given by
\[
    p(i,j,k) = V_{ij}^k = N_{ij}^k \cdot K, \quad 
    \A(i,j) = \delta_{ij} \cdot K, \quad \text{and} \quad j(i) = \delta_{0i}
    \cdot K.
\]
Braidings are described by appropriate sets of (coherent) isomorphisms
\[
    p(i,j,k) \cong p(j,i,k),
\]
associativity by isomorphisms
\[
    \bigoplus_x p(i,j,x) \otimes p(x,k,l) \cong
    \bigoplus_x p(i,x,l) \otimes p(j,k,x),
\]
and $*$-autonomy by the cyclic condition
\[
    p(i,j,Sk)=p(j,k,Si)=p(k,i,Sj)
\]
where
\[
    S:\A \ra \A, S^2 = 1
\]
is the involution of the RCFT.

Here we insist also that
\[
    p(Si,Sj,Sk) = p(i,j,k)^*,
\]
where $p(i,j,k)^*$ is the dual space of $p(i,j,k)$.

A useful way of abstracting this situation, especially for the purposes of
constructing rigid convolutions of the form $[\A,\Vect_{\text{fd}}]$, is to
replace the finite set $\A$ above by any promonoidal
$\Vect_{\text{fd}}$-category $(\A,p,j)$ with $\ob\A$ finite; then this $\A$
has a distinguished base object if $\A$ has an identity object representing
$j$.
\end{example}

\begin{remark} \quad
\begin{enumerate}[{\upshape (i)}]
\item $*$-autonomous monoidal categories (under that name) were introduced by M. Barr~\cite{1}, and then studied extensively by many authors. Their relationship to classical Frobenius structures was recognized in~\cite{4} and in~\cite{9}.

\item The $*$-autonomous structure on the extended real numbers was noted
(by the author) in the context of a lecture entitled ``$*$-autonomous
convolution'' (Australian Category Seminar, March 5, 1999), and recently
introduced anew by M.\ Grandis for other purposes.

\item The examples above admit generalization to more elaborate
promonoidal settings. Are there corresponding physical interpretations?
\end{enumerate}
\end{remark}

\subsubsection*{Acknowledgements}
I wish to thank Craig Pastro and Ross Street for general assistance with typesetting and the Australian Research Council for some financial support.



\end{document}